\documentclass[a4paper,12pt]{article}
\usepackage{amsmath,amssymb,amsthm}
\usepackage{epsfig,graphicx,verbatim}

\numberwithin{equation}{section}
\newtheorem{thm}[equation]{Theorem}

\newtheorem{lem}[equation]{Lemma}

\newcounter{mycount}
\newenvironment{romlist}{\begin{list}{\rm(\roman{mycount})}%
   {\usecounter{mycount}\labelwidth=1cm\itemsep 0pt}}{\end{list}}

\def\EE{{\mathbb E}}

\def\sM{{\mathcal M}}

\def\ZRC{Z^{\text{RC}}}
\def\s{\sigma}

\def\qq{\qquad}
\def\q{\quad}

\def\b{\beta}
\def\d{\delta}
\def\l{\lambda}
\def\g{\gamma}
\def\t{\theta}
\def\rc{random-cluster}
\def\ZZ{{\mathbb Z}}
\def\RR{{\mathbb R}}
\def\LL{{\mathbb L}}
\def\Pr{{\mathbb P}}
\def\Qr{{\mathbb Q}}
\def\la{\langle}
\def\ra{\rangle}

\def\Om{\Omega}
\def\om{\omega}

\def\La{\Lambda}
\def\oo{\infty}

\def\pc{p_{\mathrm c}}

\def\be{\begin{equation}}
\def\ee{\end{equation}}
\def\sm{\setminus}

\def\ci{\cite}
\def\ov#1{{\vec #1}}

\def\ZP{Z^{\text{P}}}
\def\pim{P}
\def\pd{\partial}
\def\lra{\leftrightarrow}
\def\es{\varnothing}
\def\symdif{\,\bigtriangleup\,}
\def\bigmid{\,\mid\,}
\def\fpq{\phi_{p,q}}
\def\fbpq{\phi_{\bp,q}}
\def\bJ{\mathbf{J}}
\def\bp{\mathbf{p}}
\def\bP{{\mathbf{P}}}
\def\bM{{\mathbf{M}}}
\def\bN{{\mathbf{N}}}
\def\bm{{\mathbf{m}}}
\def\bl{{\boldsymbol{\lambda}}}
\def\bn{{\mathbf{n}}}

\begin{document}
\title{Flows and ferromagnets\footnote{Version of 21 July 2005.}}
\author{Geoffrey Grimmett\\
Statistical Laboratory, Centre for
Mathematical Sciences,\\
University of Cambridge,\\
Wilberforce Road, Cambridge CB3 0WB, U.K.}
\maketitle

\begin{abstract}
The two-point correlation function of a Potts model
on a graph $G$ may be expressed in terms of
the flow polynomials of `Poissonian'
random graphs derived from $G$ by replacing each
edge by a Poisson-distributed number of copies of itself.
This fact extends to Potts models the so-called
random-current expansion of the Ising model.
\end{abstract}

\section{Introduction}
The Tutte polynomial and its relatives have rarely been distant from the work of Dominic Welsh. They play important roles in
matroid theory, \ci{Wel76}, computational complexity, 
\ci{Wel94,Wel99,WeMe},
and models of statistical physics, \ci{Wel93,WeMe}. They
provide the natural way to count and relate a variety of
objects defined on graphs. We show here
that they permit a representation of the two-point correlation
function of a ferromagnetic Potts model on a graph $G$ 
in terms of the flow polynomials of certain related random graphs.
This representation extends to general Potts models
the so-called random-current expansion for Ising models,
wielded with great effect in \ci{A82,ABF,AF,Si80}
and elsewhere, and it
amplifies the links 
between the Potts partition function and the Tutte polynomial
surveyed earlier by Welsh and Merino, \ci{WeMe}. 

Two key elements of the analysis of the Ising model
on a graph $G$ are the \rc\ representation and the random-current expansion. The former is valid for all Potts models (and 
more besides), but the latter has not previously been
extended beyond the Ising model.  It hinges on an expansion of the
partition function in terms of $0/1$-vectors indexed by edges 
and such that, for every vertex $v$,
the sum of the values over edges incident to $v$ is even. Such 
a  vector
may be recognised as a `mod-2 flow'. It turns out that the $q$-state Potts partition function corresponds similarly
to counts of `mod-$q$ flows' on a graph derived
from $G$ in the following way. Let $\lambda>0$, and replace every
 edge $e$ of $G$ by $P(e)$ parallel
edges, where the $P(e)$ are independent Poisson-distributed
random variables with parameter $\lambda$. 
The quantity of interest is the mean
number of non-zero mod-$q$ flows on the resulting
random graph.

There is a powerful method of `path-manipulation' by which
many important results have been 
proved for the Ising model. This method
has a simple form when set in the context of a Poissonian random
graph, and we illustrate this in Section 5 with a version
of the `switching lemma' of \ci{A82}.

A short tour of graph polynomials appears in Section 2. In Section
3 is introduced the Potts and \rc\ models, and the main result is
proved in Section 4. Applications to the Ising
model are summarized in Section 5. The principal open
area is to extend the random-current analysis to Potts models
with general $q$.

\section{Graph polynomials}
Let $G=(V,E)$ be a finite graph, possibly containing multiple
edges and loops. The Whitney and Tutte polynomials
of $G$ are well known to 
graph theorists, and we begin with a reminder of 
their definitions.
The (Whitney) {\it rank-generating function\/} of $G$
was introduced in \ci{W32} and is given by
\begin{equation}
W_G(u,v)=\sum_{E'\subseteq E} u^{r(G')} v^{c(G')},\qq u,v\in\RR,
\label{old2.74}
\end{equation}
where $r(G')=|V|-k(G')$ is the {\it rank\/} of the subgraph $G'=(V,E')$,
and $c(G')=|E'|-|V|+k(G')$ is its {\it co-rank\/}. Here, $k(G')$
denotes the number of components of $G'$. 
Note that
\begin{equation}
W_G(u,v) = (u/v)^{|V|} \sum_{E'\subseteq E} v^{|E'|} (v/u)^{k(E')},
\qq u,v \ne 0.
\label{8.rank2}
\end{equation}

The rank-generating function has various useful properties, and 
it occurs in several contexts in graph theory, see \ci{Big,Tut}. 
The Tutte (or dichromatic) polynomial of $G$ was introduced 
independently in
\ci{Tut47,Tut}, and may be expressed as
\begin{equation}
T_G(u,v)= (u-1)^{|V|-1} W_G\bigl( (u-1)^{-1}, v-1\bigr).
\label{old2.75}
\end{equation}
This also is a function of two variables.
For suitable values of these variables,
it provides counts  of
colourings, forests, and flows, and of other combinatorial
quantities. The principal purpose of the current paper is to explore
the use of the Whitney/Tutte polynomial in the study of the correlation
functions of the Potts model, and to this end we define next the 
{\it flow polynomial\/} of $G$.

We turn $G$ into a oriented graph by allocating
an orientation to each edge $e\in E$, 
and we denote the resulting
digraph by $\ov G=(V, \ov E)$. If the edge $e=\la u,v\ra\in E$ is 
oriented from $u$ to $v$, 
we say that $f$ {\it leaves\/} $u$ and 
{\it arrives\/} at $v$.
It will turn out that the choices
of orientations are immaterial to the principal conclusions that follow. Let $q\in\{2,3,\dots\}$.
A function $f:\ov E\to\{0,1,2,\dots,q-1\}$ is called
a {\it mod-$q$ flow\/} on $\ov G$ if
$$
\sum_{\text{\scriptsize$\ov e\in\ov E:$}\atop \text{\scriptsize$\ov e$ leaves $v$}} f(\ov e) -
\sum_{\text{\scriptsize$\ov e\in\ov E:$}\atop \text{\scriptsize$\ov e$ arrives at $v$}} f(\ov e) =
0\qq\text{modulo $q$,\q for all $v\in V$},
\label{8.flowdef}
$$
which is to say that flow is conserved at every vertex.
A mod-$q$ flow $f$ is called {\it non-zero\/} if $f(\ov e)\ne 0$ for all
$\ov e\in \ov E$. 
Let $C_G(q)$ be the number of non-zero mod-$q$ flows on $\ov G$.
It is fundamental that 
the quantity $C_G(q)$ does not depend on the orientations
of the edges of $G$, and the proof may be found in
\ci{Tut}.
The function $C_G(q)$, viewed as a function of $q$, is 
called the {\it flow polynomial\/} of $G$.

The flow polynomial may be obtained
as an evaluation of the Whitney/Tutte polynomial with
two particular parameter values, as follows:
\begin{align}
C_G(q)&=(-1)^{|E|} W_G(-1,-q)
\nonumber\\
&= (-1)^{|E|-|V|+1}T_G(0,1-q),
\qq q\in\{2,3,\dots\}. \label{8.flowpodef} 
\end{align}
See \ci{Big,Tut}. We shall later write $C(G;q)$ for
$C_G(q)$, and similarly for other polynomials when the
notational need arises.

\section{Potts and \rc\ models}\label{secPoIs}
Amongst models for ferromagnetism, 
the Potts model is one of the most studied.
It has two principal parameters, the `inverse temperature'
$\b\in(0,\oo)$
and the number $q\in\{2,3,\dots\}$ of local states. 
When $q=2$, the
Potts model becomes the Ising model. 
Let $G=(V,E)$ be a finite graph which for simplicity
we assume to have no loops.
It is convenient to allow a separate parameter for each edge of $G$,
and thus we let $\bJ=(J_e: e\in E)$ be a  vector of non-negative 
numbers, and we set 
\begin{equation}
p_e=1-e^{-\beta J_e q}, \qq e \in E.
\label{8.pj}
\end{equation}
The configuration space of the $q$-state Potts model on $G$ is
the set $\Sigma = \{1,2,\dots,q\}^V$. The Potts measure on 
$\Sigma$ is given by
\begin{equation}
\pi_{\beta \bJ,q} (\s) = \frac 1{\ZP} \exp\biggl\{\sum_{e \in E} \beta J_e(q\delta_e(\s)-1)\biggr\},
\qq \s\in \Sigma,
\label{8.Pottsdef}
\end{equation}
where, for $e=\la x,y\ra \in E$,
$$
\d_e(\s) = \d_{\s_x,\s_y} = \begin{cases}
 1 &\text{if } \s_x=\s_y,\\
0 &\text{otherwise,}
\end{cases}
$$
and $\ZP=\ZP_G$ is the partition function
\begin{equation}
\ZP= \sum_{\s\in\Sigma} \exp\biggl\{\sum_{e \in E} 
\beta J_e(q\delta_e(\s)-1)\biggr\}.
\label{8.pottspart}
\end{equation}
Since $\b J_e\ge 0$, the Potts measure $\pi_{\b\bJ,q}$
allocates greater probability to configurations for which
$\d_e(\s)=1$ for a larger set of edges $e$.
That is, it prefers
configurations in which
many neighbour-pairs have the same state, and in this regard the
model is termed `ferromagnetic'.

A central quantity is
the `two-point correlation function' given 
by
\begin{equation}
\tau_{\beta \bJ,q}(x,y) = \pi_{\beta \bJ,q}(\s_x=\s_y)-\frac 1q,
\qq x,y \in V.
\label{8.taudef}
\end{equation}
We shall work here with $q\tau_{\beta \bJ,q}(x,y)=\pi_{\beta \bJ,q}
(q\d_{\s_x,\s_y}-1)$ and, for ease of notation in the following, we write
\begin{equation}
\s(x,y)=q\tau_{\beta \bJ,q}(x,y),\qq x,y\in V,
\label{8.sigdef}
\end{equation}
thereby suppressing reference to the parameters $\beta \bJ$ and $q$.

Two of the most successful ways of studying the Ising/Potts
models are the so-called `\rc\ model' and the `random-current
expansion'. 
We define next the \rc\ model, and we explain its relevance
to the Potts model. The random-current expansion
for the Ising model will be
reviewed in Section \ref{secrcurr}.

In the (bond) percolation model on $G$, 
each edge is declared at random
to be either `open' or `closed'. 
An edge is declared `open' with some given
probability $p$, and closed otherwise, and different edges 
are allocated independent states. The percolation model
is basic to the study of disordered media, particularly
when the underlying graph is part of a `crystalline'
lattice such as the $d$-dimensional cubic lattice $\LL^d$.
See \ci{G99} for a full account. When $G$ is a complete graph, the percolation model is usually called an `Erd\H os--R\'enyi
random graph', see \ci{JLR}.

The \rc\ measure on $G$ is obtained through a perturbation
of the percolation measure, as follows. Let $\bp
=(p_e: e\in E) \in[0,1]^E$
and $q\in(0,\oo)$. The configuration space is $\Om=\{0,1\}^E$.
For $\om\in\Om$ and $e\in E$, we say that $e$ is $\om$-{\it open\/}
(or, simply, open) if $\om(e)=1$,
and $\om$-{\it closed\/} otherwise.
The \rc\ probability measure on $\Om$ is defined by
$$
\fbpq(\om) = \frac1{\ZRC} 
\biggl\{ \prod_{e\in E} p_e^{\om(e)} (1-p_e)^{1-\om(e)}
                           \biggr\} q^{k(\om)}, \qq \om \in \Om, 
\label{old2.53}
$$
where $k(\om)$ denotes the number of $\om$-open
components on the vertex-set $V$, and
$\ZRC=\ZRC_G$ is the appropriate normalizing factor. 
We sometimes write $\phi_{G,\bp,q}$ when the role of $G$ is to be emphasized. 

It is common to take $p_e=p$ for all $e\in E$, in
which case we write $\fpq$ for $\fbpq$. The special case
$q=1$, $\bp=p$ is evidently the percolation measure
with parameter $p$, in which case we write
$\phi_p=\phi_{p,1}$. It turns out that the \rc\ model
with $q\in\{2,3,\dots\}$ corresponds in a certain way to the
Potts model on $G$ with $q$ local states and with
$\b\bJ$ satisfying \eqref{8.pj}. Specifically,
the two-point correlation function of the latter is 
(up to a harmless factor)
equal to the connection probability of the former,
\begin{equation}
\tau_{\b\bJ,q}(x,y) = (1-q^{-1})\fbpq(x\lra y),
\qq x,y\in V,
\label{corrconn}
\end{equation}
where $x\lra y$ means that there exists a path of
open edges from $x$ to $y$.
The \rc\ model was introduced by Fortuin and Kasteleyn
around 1970, and has been reviewed recently
in \ci{G95,G02,G-RC}.

The \rc\ partition function $\ZRC_G$
is given by
$$
\ZRC_G(p,q)=\sum_{\om\in\Om} p^{|\eta(\om)|}(1-p)^{|E\setminus\eta(\om)|}
q^{k(\om)},
\label{old2.73}
$$
and is easily seen by \eqref{8.rank2} to satisfy
$$
\ZRC_G(p,q) = q^{|V|} (1-p)^{|E|} W_G
\left(\frac p{q(1-p)}, \frac p{1-p}\right) ,
\qq p\ne 1,
\label{old2.76}
$$
a relationship which provides a link with other classical 
graph-theoretic quantities. See \ci{Big77,Big,F72a,Sok05,WeMe}.

\section{Potts correlations and flow counts}
It is shown in this section
that the  Potts correlation functions \eqref{8.taudef}
may be expressed in terms of flow polynomials associated
with a certain `Poissonian' random graph derived from $G$
by replacing each edge by a random number of copies.
This extends to general $q$ the random-current expansion of the Ising model described in Section \ref{secrcurr}.

For
any vector $\bm=(m(e): e\in E)$ of non-negative integers,
let $G_\bm=(V,E_\bm)$ be the graph with vertex
set $V$ and, for each $e\in E$, with exactly $m(e)$ edges
in parallel joining the endvertices of the edge $e$;
the original edge $e$ is itself removed.
Note that
\begin{equation}
|E_\bm| = \sum_{e\in E} m(e).
\label{8.sumem}
\end{equation}
Let $\bl=(\l_e:e\in E)\in[0,\oo)^E$. Let $\bP=
(\pim(e):e\in E)$ be a family of independent random variables
such that $\pim(e)$ has the Poisson 
distribution with parameter $\l_e$.
The random graph $G_\bP=(V,E_\bP)$ is
 called a {\it Poisson graph with intensity $\bl$}. Let
$\Pr_\bl$ and $\EE_\bl$ denote the corresponding probability measure
and expectation operator.

For $x,y\in V$, $x\ne y$,
we denote by $G_\bP^{x,y}$ the graph obtained from $G_\bP$ by adding an edge
with endvertices $x$, $y$. If $x$ and $y$ are already 
adjacent in $G_\bP$, we add exactly
one further edge between them. 
Potts-correlations and flows are related by the following theorem.
The function $\s(x,y)$ is given in \eqref{8.sigdef}.

\begin{thm}\label{8.flowcorr} 
Let $q \in\{2,3,\dots\}$ and $\l_e=\beta J_e$. Then
\begin{equation}
\s(x,y) = \frac{\EE_\bl(C(G_\bP^{x,y};q))}
{\EE_\bl(C(G_\bP;q))},
\qq x,y\in V.
\label{8.flowcorr2}
\end{equation}
\end{thm}

This formula is particularly striking when $q=2$, since
non-zero mod-2 flows necessarily take only the value 1. 
A finite graph $H=(W,F)$ is called {\it even\/} if the degree of every
vertex $w$ is even.
It is trivial that $C_H(2)=1$ 
if $H$ is even, and $C_H(2)=0$ otherwise.
By \eqref{8.flowcorr2},
for any graph $G$, 
\begin{equation}
\s(x,y)=\frac{\Pr_\bl(G^{x,y}_\bP\text{ is even})}
{\Pr_\bl(G_\bP\text{ is even})}.
\label{8.flowcorr3}
\end{equation}
Such observations are at the heart of the
random-current expansion for Ising models.
See Section \ref{secrcurr}.

Theorem \ref{8.flowcorr}
may be extended via \eqref{corrconn} to the \rc\ model.
Assume for simplicity that every edge has the same
parameter $p$. The proof of the following is
easily derived from Theorem \ref{8.flowcorr}, and may be found in 
\ci{G-RC}. It is obtained by expressing the flow polynomial
in terms of the Tutte polynomial $T$, and allowing $q$ to 
vary continuously.

\begin{thm}\label{8.flowconn} 
Let $p\in[0,1)$ and $q \in(0,\oo)$. Let $\l_e=\l$
for all $e\in E$, where
$p=1-e^{-\lambda  q}$. 
\begin{romlist}
\item 
For $x,y\in V$,
\begin{equation}
(q-1)\phi_{G,p,q}(x\lra y) = 
\frac{\EE_\l\bigl((-1)^{1+|E_\bP|} T(G_\bP^{x,y};0,1-q)\bigr)}
{\EE_\l\bigl((-1)^{|E_\bP|}T(G_\bP;0,1-q)\bigr)}.
\label{8.flowconn2}
\end{equation}
\item
For $q\in\{2,3,\dots\}$,
\begin{equation}
\phi_{p}(q^{k(\om)})=
(1-p)^{|E|(q-2)/q} q^{|V|} \EE_\l(C(G_\bP;q)).
\label{8.compflow}
\end{equation}
\end{romlist}
\end{thm}

When $q=2$, \eqref{8.compflow} reduces to the curiosity
\begin{equation}
\phi_{p}(2^{k(\om)}) = 2^{|V|} \Pr_\l(G_\bP \text{ is even}).
\label{8.compflow2}
\end{equation}
This may be simplified further. Let $\zeta(e)=\pim(e)$ modulo 2. It is easily seen that $G_\bP$ is an even graph
if and only if $G_{\zeta}$ is even, and that the $\zeta(e)$, $e\in E$,
are independent Bernoulli variables with
$$
\Pr_\l(\zeta(e)=1) = \tfrac12(1-e^{-2\l}) = \tfrac12 p.
$$
Equation \eqref{8.compflow} may therefore be written as
\begin{equation}
\phi_{p}(2^{k(\om)}) = 2^{|V|} 
\phi_{p/2}(\text{the open graph on $V$ is even}).
\label{8.star}
\end{equation}

\begin{proof}[Proof of Theorem \ref{8.flowcorr}]
Since the parameter $\beta$ appears together
with the multiplicative factor $J_e$, 
we may without loss of generality
take $\beta=1$. 
We begin with a calculation involving the Potts 
partition function $\ZP$ of
\eqref{8.pottspart}. Let $\ZZ_+=\{0,1,2,\dots\}$
and let $\bm=(m_e:e\in E)\in\ZZ_+^E$.
By a Taylor expansion in the variables $J_e$,
\begin{align}
\exp\biggl\{-\sum_{e\in E} J_e\biggr\} \ZP &=
  \sum_{m \in \ZZ_+^E} \biggl( \prod_{e\in E} \frac{J_e^{m_e}}{m_e!}
  e^{-J_e}\biggr) \pd^\bm\ZP\Big|_{\bJ=0} \nonumber
\\
&= \EE_\bl\left(\pd^\bP\ZP\Big|_{\bJ=0}\right)
\label{8.taylor}
\end{align}
where
$$
\pd^\bm\ZP = \left(\prod_{e\in E} \frac{\pd^{m_e}}{\pd J_e^{m_e}}\right)\ZP,
\qq \bm\in\ZZ_+^E.
$$
By \eqref{8.pottspart} with $\b=1$, 
\begin{align}
\pd^\bm\ZP\Big|_{\bJ=0} &=\sum_{\s\in\Sigma}\,\prod_{e\in E} (q\d_e(\s)-1)^{m_e} \nonumber\\
&=\sum_{\s\in\Sigma}\, \prod_{e\in E_\bm}  (q\d_e(\s)-1)\nonumber\\
&= \sum_{\s\in\Sigma}\, \prod_{e\in E_\bm}\, \sum_{n_e\in\{0,1\}} 
  [-\d_{n_e,0} + \d_{n_e,1} q \d_e(\s)]\nonumber\\
&= \sum_{\bn\in\{0,1\}^{E_\bm}}\, \sum_{\s\in\Sigma} (-1)^{|\{e:\,n_e=0\}|}
 q^{|\{e:\,n_e=1\}|} \biggl(\prod_{e\in E_m} \d_e(\s)^{n_e}\biggr)
\nonumber\\
&= \sum_{\bn\in\{0,1\}^{E_\bm}}  (-1)^{|\{e:\,n_e=0\}|} q^{|\{e:\,n_e=1\}|}
  q^{k(m,n)},\label{8.zpart} 
\end{align}
where $k(\bm,\bn)$ is the number of connected components of
the graph obtained from $G_\bm$ after deletion of all edges $e$ with $n_e=0$.
By \eqref{8.rank2}--\eqref{8.flowpodef},
\begin{align}
\pd^\bm\ZP\Big|_{\bJ=0} &= (-1)^{|E_\bm|} 
\sum_{\bn\in\{0,1\}^{E_\bm}} (-q)^{|\{e:\,n_e=1\}|}
q^{k(\bm,\bn)}\nonumber  \\
&= (-1)^{|E_\bm|} q^{|V|} W_{G_\bm}(-1,-q)\\
&= q^{|V|} C(G_\bm;q). \label{8.zpart2} 
\end{align}
Combining \eqref{8.taylor}--\eqref{8.zpart2}, 
\begin{equation}
\exp\biggl\{-\sum_{e\in E} J_e\biggr\} \ZP = q^{|V|} \EE_\bl(C(G_\bP;q)).
\label{8.partflow}
\end{equation}

Let $x,y\in V$. We define the unordered pair $f=(x,y)$,
and write $\d_f(\s)=\d_{\s_x,\s_y}$ for $\s\in\Sigma$.
We have that
\begin{align}
\s(x,y) &= \pi_{\beta\bJ,q}(q\d_f(\s)-1)
\nonumber \\
&=\frac 1{\ZP} \sum_{\s\in\Sigma}(q\d_f(\s)-1) \exp\biggl\{\sum_{e \in E}
\beta J_e(q\delta_e(\s)-1)\biggr\}.
\label{8.partflow3} \end{align}
By an analysis parallel to 
\eqref{8.taylor}--\eqref{8.partflow}, 
\begin{align}
\exp\biggl\{-\sum_{e\in E} J_e\biggr\} 
&\sum_{\s\in\Sigma}(q\d_f(\s)-1) \exp\biggl\{\sum_{e \in E}
\beta J_e(q\delta_e(\s)-1)\biggr\} \label{8.partflow2}\\
&\hskip.5cm = q^{|V|} \EE_\bl(C(G_\bP^{x,y};q)),
\nonumber
\end{align}
and \eqref{8.flowcorr2} follows by \eqref{8.partflow} and 
\eqref{8.partflow3}.
\end{proof}

\section{Random-current expansion of the Ising model}\label{secrcurr}
Unlike the situation with the Potts model, 
there is a fairly complete
analysis of the Ising model. A principal part
in this analysis is played by Theorem \ref{8.flowcorr} with $q=2$,
under the heading `random-current expansion'. This has permitted
proofs amongst other things of the exponential decay
of correlations in the low-$\b$ regime on the cubic lattice
$\LL^d$ with $d\ge 2$. See \ci{A82,ABF,AF}. It has not so far been
possible to extend this work to general Potts models, but 
Theorem \ref{8.flowcorr}
could play a part in such an extension.

Let $G=(V,E)$ be a finite graph without loops 
as before, and set $q=2$.
We restrict ourselves here to the Ising model with
$J_e=J$ for all $e\in E$, and we write $\l=\b J$.
By Theorem \ref{8.flowcorr},
\begin{equation}
\s(x,y)=2\tau_{\l,2}(x,y)=\frac{\Pr_\l(\text{$G_\bP^{x,y}$ is even})}
{\Pr_\l(\text{$G_\bP$ is even})},
\qq 0\le \l<\oo.
\label{8.isingeven}
\end{equation}
The value of such a representation will become clear 
during the following discussion, 
which is based on material in
\ci{A82,Li80,Si80}. 
In advance of this, we make a remark concerning
\eqref{8.isingeven}. In deciding whether $G_\bP$ or $G_\bP^{x,y}$
is an even graph, we need only know the numbers $\pim(e)$ 
when reduced modulo 2.
That is, we can work with $\zeta\in\Omega=\{0,1\}^E$ given by
$\zeta(e)=\pim(e)$ mod $2$. Since $\pim(e)$ has the Poisson distribution 
with parameter $\l$, $\zeta(e)$ has the Bernoulli distribution with parameter
$$
p'=\Pr_\l(\pim(e)\text{  is odd}) = \tfrac12(1-e^{-2\l}).
$$
We obtain thus from \eqref{8.isingeven} that
$$
\s(x,y)=\frac{\phi_{p'}(\pd\zeta=\{x,y\})}{\phi_{p'}(\pd\zeta=\es)},
\label{8.isingperc}
$$
where $\phi_{p'}$ denotes product measure on 
$\Om$ with density $p'$, and
$$
\pd\zeta = \biggl\{v\in V: 
\sum_{e:\,e\sim v} \zeta(e)\text{ is odd}\biggr\},
\qq\zeta\in\Om,
$$
where the sum is over all edges $e$ incident to $v$.
We refer to members of $\pd\zeta$ as `sources' of the 
configuration $\zeta$.

Let $\bM=(M_e:e\in E)$
be a sequence of disjoint finite sets (possibly empty)
indexed by $E$, and let $m_e=|M_e|$.
As noted in the last section,
the vector $\bM$ may be used to construct a 
multigraph $G_\bm=(V,E_\bm)$ in which
each $e\in E$ is replaced by $m_e$ edges in parallel; we may
take $M_e$ to be the set of such edges. For $x,y\in V$, we write
`$x \lra y$ in $\bm$' if $x$ and $y$ 
lie in the same component of $G_\bm$. 
We define the set $\pd \bM$ of {\it sources\/} of $\bM$ by
\begin{equation}
\pd \bM= \biggl\{v\in V: \sum_{e:\,e\sim v} m_e \text{ is odd}\biggr\}.
\label{8.source1}
\end{equation}
Thus, for example,
$G_\bm$ is even if and only if 
$\pd \bM$ is empty. From the vector $\bM$ we
construct a vector $\bN=(N_e: e\in E)$ by deleting each member of each $M_e$
with probability $\frac12$, independently of all other elements. That is,
we let $B_i$, $i\in\bigcup_e M_e$, be independent Bernoulli random variables
with parameter $\frac12$, and we set
$$
N_e = \{i \in M_e: B_i=1\}, \qq e \in E.
$$
We write $\Pr^\bM$ for the appropriate probability measure.
The following lemma is pivotal for the computations which follow.

\begin{lem}\label{8.switching}
Let
$\bM$ and $\bm$ be as above. If $x,y\in V$
are such that $x\ne y$ and $x\lra y$ in $\bm$ then, for $A\subseteq V$,
$$
\Pr^\bM\bigl(\pd \bN= \{x,y\},\, \pd(\bM\sm \bN)= A\bigr)
= \Pr^\bM\bigl(\pd \bN= \es ,\, \pd(\bM\sm \bN)=A\symdif\{x,y\}  \bigr).
$$
\end{lem}

\begin{proof}
Take $M_e$ to be the set of edges
of $G_\bm$ parallel to $e$, and assume that $x\lra y$ in $\bm$.
Fix $A\subseteq V$.
Let $\sM$ be the set of all vectors $\bn=(n_e:e\in E)$ with
$n_e\subseteq M_e$ for all $e$.
Let $\bp$ be a fixed path of $G_\bm$ with endpoints  $x$, $y$, and consider
the map $\rho:\sM\to\sM$ given by
$$
\rho(\bn) = \bn\symdif \bp, \qq \bn\in \sM. 
$$
The map $\rho$ is  one--one, and maps 
$\{\bn\in\sM: \pd \bn= \{x,y\} ,\ \pd(\bM\sm \bn)=A \}$
to $\{\bn'\in\sM: \pd \bn'= \es ,\ \pd(\bM\sm \bn')=A\symdif \{x,y\}  \}$.
Each member of $\sM$ is equiprobable under $\Pr^\bM$, and the claim follows.
\end{proof}

Let $\l \in[0,\oo)$, and recall
from the last section the definition of a Poisson graph
with parameter $\l$.
The following is a fairly immediate corollary of the last theorem.
Let $\bM=(M_e:e\in E)$ and $\bM'=(M_e':e\in E)$ be vectors of disjoint
finite sets satisfying $M_e\cap M_f'=\es$ for
all $e,f \in E$, and suppose that the random variables
$m_e=|M_e|$, $m_e'=|M_e'|$, $e\in E$ are independent and such that,
for each $e\in E$,
$m_e$ and $m_e'$ have the Poisson distribution with parameter $\l$.
We write $\bM\cup \bM'=(M_e\cup M_e': e\in E)$, and
$\Pr$ for the appropriate probability measure.
The following lemma is a simplification
of the so-called switching lemma of \ci{A82}.

\begin{lem}\label{8.switch2}
If $x,y\in V$
are such that $x\ne y$ and $x\lra y$ in $\bm+\bm'$ then, for $A\subseteq V$,
\begin{multline}
\Pr\bigl(\pd \bM=\{x,y\},\, \pd \bM'=A\bigmid \bM\cup \bM'\bigr)\\
= \Pr\bigl(\pd \bM=\es,\, \pd \bM'=A\symdif\{x,y\}\bigmid \bM\cup \bM'\bigr)
\qq\Pr\text{\rm -a.s.}
\end{multline}
\end{lem}

\begin{proof}
Conditional on the sets $M_e\cup M'_e$ for $e\in E$,
the sets $M_e$ are selected by the 
independent removal of each element with probability $\frac12$.
The claim now follows from Lemma \ref{8.switching}.
\end{proof}

We present two applications of Lemma
\ref{8.switch2} to the Ising model,
as in \ci{A82}. For $\bm=(m_e:e\in E)\in \ZZ_+^E$, let
\begin{equation}
\pd \bm=
 \biggl\{v\in V: \sum_{e:\,e\sim v} m_e \text{ is odd}\biggr\},
\label{8.source2}
\end{equation}
as in \eqref{8.source1}. We write as before
$$
\s(x,y)=2\tau_{\l,2}(x,y)=\pi_{\l,2}(2\d_{\s_x,\s_y}-1),\qq x,y\in V,
$$
thereby suppressing reference to $\l$.
By \eqref{8.isingeven},
\begin{equation}
\s(x,y)= \frac{\Pr_\l(\pd\bP=\{x,y\})}{\Pr_\l(\pd\bP=\es)}.
\label{8.tausource}
\end{equation}
Let $\Qr_A$ denote the law of $\bP$ {\it
conditional\/} on the event $\{\pd\bP=A\}$, 
$$
\Qr_A(E) = \Pr_\l(\bP\in E\mid \pd\bP=A).
$$
We have need of two independent copies $\bP_1$, $\bP_2$ of $\bP$ with
potentially
different conditionings, and thus we write $\Qr_{A;B}=\Qr_A\times \Qr_B$.

\begin{lem}\label{8.appls}
Let $x,y,z\in V$ be distinct vertices. Then\/{\rm:}
\begin{romlist}
\item  
$\s(x,y)^2 = \Qr_{\es;\es}(x\lra y\text{ \rm in } \bP_1+\bP_2)$,
\item
$\s(x,y)\s(y,z)
= \s(x,z)\Qr_{\{x,z\};\es}(x\lra y\text{ \rm in } \bP_1+\bP_2)$.
\end{romlist}
\end{lem}

\begin{proof}
(i) By \eqref{8.tausource} and Theorem \ref{8.switch2},
\begin{align*}
\s(x,y)^2 &= \frac{\Pr_\l\times\Pr_\l(\pd\bP_1=\{x,y\},\,
                           \pd\bP_2=\{x,y\})}
			   {\Pr_\l(\pd\bP=\es)^2}\\
&=
\frac{\Pr_\l\times\Pr_\l
(\pd\bP_1=\{x,y\},\,\pd\bP_2=\{x,y\},\, x\lra y \text{ in } \bP_1+\bP_2)}
{\Pr_\l(\pd\bP=\es)^2}\\
&= 
\frac{\Pr_\l\times\Pr_\l
(\pd\bP_1=\pd\bP_2=\es,\, x\lra y \text{ in } \bP_1+\bP_2)}
{\Pr_\l(\pd\bP=\es)^2}\\
&=\Qr_{\es;\es}(x\lra y \text{ in } \bP_1+\bP_2).
\end{align*}

\noindent
(ii) Similarly,
\begin{align*}
&
\s(x,y)\s(y,z)\\
&\q= \frac{\Pr_\l\times\Pr_\l\bigl(\pd\bP_1=\{x,y\},\, \pd\bP_2=\{y,z\}\bigr)}
{\Pr_\l(\pd\bP=\es)^2}\\
&\q=
\frac{\Pr_\l\times\Pr_\l\bigl(\pd\bP_1=\es,\, \pd\bP_2=\{x,z\},\,
x\lra y \text{ in } \bP_1+\bP_2\bigr)}
{\Pr_\l(\pd\bP=\es)^2}\\
&\q= \frac{\Pr_\l(\pd\bP_2=\{x,z\})}{\Pr_\l(\pd\bP=\es)}
\cdot 
\Pr_\l\times\Pr_\l
  \bigl(x\lra y \text{ in } \bP_1+\bP_2\bigmid \pd\bP_1=\es,\ 
  \pd\bP_2=\{x,z\}\bigr)\\
&\q=
\s(x,z)\Qr_{\{x,z\};\es}(x\lra y\text{ in } \bP_1+\bP_2),
\end{align*}
and the proof is complete.
\end{proof}

Theorem \eqref{8.appls}(ii) leads to an important correlation inequality
known as the `Simon inequality', \ci{Si80}. 
Let $x,z\in V$ be distinct vertices.
A subset $W\subseteq V$ is said to {\it separate\/} $x$ and $z$ if 
$x,z\notin W$ and every path from $x$ to $z$ contains some vertex of $W$.

\begin{thm}\label{8.simon}
Let $x,z\in V$ be distinct vertices, and let $W$ separate
$x$ and $z$. Then
$$
\s(x,z)\le \sum_{y\in W} \s(x,y)\s(y,z).
$$
\end{thm}

\begin{proof}
By Theorem \ref{8.appls}(ii),
\begin{align*}
\sum_{y\in W} \frac{\s(x,y)\s(y,z)}{\s(x,z)}
&= \sum_{y\in W} \Qr_{\{x,z\};\es}(x\lra y \text{ in } \bP_1+\bP_2)\\
&= \Qr_{\{x,z\};\es}\bigl(\bigl|\{y\in W: x\lra y\text{ in } \bP_1+\bP_2\}
\bigr|\bigr).
\end{align*}
Assume that the event  $\pd\bP_1=\{x,z\}$ occurs.
On this event,  $x\lra z$ in $\bP_1+\bP_2$.
Since $W$ separates $x$ and $z$, 
the set $\{y\in W: x\lra y\text{ in } \bP_1+\bP_2\}$ is 
non-empty on this event. Thus its (conditional) mean 
size is at least one under $\Qr_{\{x,z\};\es}$,
and the claim
follows.
\end{proof}

The Ising model on the graph
$G=(V,E)$ corresponds as described
in Section \ref{secPoIs} to a \rc\ measure $\phi_{p,q}$ with $q=2$.
By \eqref{corrconn},
$$
\s(x,y)=2 \tau_{\l,2}(x,y)=\phi_{p,q}(x\lra y),
$$
where $p=1-e^{-\l q}$ and $q=2$. The Simon inequality
may be written in the form
$$
\phi_{p,q}(x\lra z)\le\sum_{y\in W}\phi_{p,q}(x\lra y)\phi_{p,q}(y\lra z)
\label{8.simon2}
$$
whenever $W$ separates $x$ and $z$. It is a curious fact that this inequality
holds also when $q=1$, see \ci{G99,H57a}. 
One may conjecture that it holds for any $q\in[1,2]$.

Let $d\ge 2$.
The \rc\ measure $\fpq$ on the cubic lattice $\LL^d$ may be
obtained as a weak limit (with so-called free boundary
conditions) of the \rc\ measure 
on finite boxes $\La$, as $\La\uparrow\ZZ^d$. The {\it
percolation probability\/} is the function $\t$ given by
$$
\t(p,q)=\fpq(0\lra \oo),
$$
the probability that the origin is the endpoint of an infinite open
path. The critical point is defined as
$$
\pc(q)=\sup\{p:\t(p,q)=0\}.
$$
Let $\|\cdot\|$ be a norm on $\ZZ^d$.
It has been conjectured that, for $p<\pc(q)$,
there exists $\g=\g(p,q)\in(0,\oo)$ such that
$\fpq(0\lra x) \le e^{-\|x\|\g}$ for all
$x \in \ZZ^d$. This has been proved when $q=1,2$ and $q$ is
sufficiently large. The Simon inequality implies the following necessary and sufficient condition for exponential decay when $q=2$.

\begin{thm}\label{8.expdecay}
Let $q=2$ and assume that $p$ is such that
$$
\sum_{x\in\ZZ^d} \fpq(0\lra x) < \oo.
\label{8.summ}
$$
There exists $\gamma=\gamma(p,q) \in (0,\oo)$ such that
$$
\fpq(0\lra z) \le e^{-\|z\|\gamma}, \qq z\in \ZZ^d.
$$
\end{thm}

The proof follows standard lines, and may be found 
in \ci{G-RC, Si80} together
with proofs of the following facts.
There is an important extension of the Simon inequality
due to Lieb, \ci{Li80}. This also may be proved via
the flow representation of Theorem \ref{8.flowcorr}.
The Lieb inequality has an important consequence for the
nature of the phase transition of the Ising model,
namely the `vanishing of the mass gap'.

Let $q=2$ and write
$$
\psi(p,q) = \lim_{n\to\oo}\left\{
-\frac 1n \log \fpq(0\lra \pd\La_n)\right\},
$$
where $\La_n=[-n,n]^d$ and $\pd\La_n=\La_n\sm\La_{n-1}$.
Note that
$\psi(p,q)$ is non-increasing in $p$, and $\psi(p,q)=0$ if $p>\pc(q)$.
One of the characteristics of a first-order phase transition
is the (strict) exponential decay of connectivity probabilities
{\it at\/} the critical point. The quantity $\psi(\pc(q),q)$
is sometimes termed the {\it mass gap\/}.

\begin{thm}\label{8.expdeclieb}
Let $d\ge 2$ and
$q=2$. Then $\psi(p,q)$ decreases to $0$ as $p\uparrow
\pc(q)$. In particular $\psi(\pc(q),q)=0$, that is,
the mass gap equals $0$.
\end{thm}

This was proved in \ci{Li80}, see also \ci{G-RC}.
The corresponding statement is known to be false for $d\ge 2$ and
$q>Q(d)$ for some sufficiently large $Q(d)$. See \ci{G-RC,LMMRS}.
In further use of the random-current expansion (with $q=2$), it
has been proved that $\psi(p,q)>0$ whenever $p<\pc(q)$.
See \ci{ABF,AF,BA91} for more details
of the Ising phase transition.

\section{Acknowledgement}
The author acknowledges the hospitality of Cornell University, 
where Theorem \ref{8.flowcorr} was proved in 1991.

\bibliography{dwgrim}
\bibliographystyle{plain}

\end{document}